\documentclass[11pt]{article}

\usepackage{amsfonts,amsmath,amssymb,amsthm}

\newtheorem{thm}{Theorem}[section]

\newtheorem{prp}[thm]{Proposition}
\newtheorem{cor}[thm]{Corollary}
\theoremstyle{definition}
\newtheorem{dfn}[thm]{Definition}
\theoremstyle{remark}
\newtheorem{rem}[thm]{Remark}

\DeclareMathOperator{\Ric}{Ric}

\begin{document}

\pagenumbering{arabic}

\begin{center}
\Large{{\bf SPECIAL STANDARD STATIC SPACE-TIMES
}  }
\end{center}

\vspace{0.5 cm}

\begin{center}
\large{Fernando Dobarro \footnote{Author partially supported by
funds of the National Group 'Analisi Reale' of the Italian
Ministry of University and Scientific Research at the University
of Trieste.} \qquad B\"ulent \"Unal }
\end{center}

\vspace{0.5 cm}

\begin{abstract}
Essentially, some conditions for the Riemannian factor and the
warping function of a standard static space-time are obtained in
order to guarantee that no nontrivial warping function on the
Riemannian factor can make the standard static space-time
Einstein.
\end{abstract}

\vspace{0.5 cm}

\renewcommand{\thepage}{\arabic{page}}
\setcounter{page}{1}

\section{Introduction}

In order to obtain general solutions to Einstein's field
equations, Lorentzian warped product manifolds were introduced
in general relativity (see \cite{BEE,ON}). Generalized
Robertson-Walker space-times and standard static space-times are
two well known important examples. The former are clearly a
generalization of Robertson-Walker space-times and the latter a
generalization of the Einstein static universe. In this paper,
we basically focus on properties of the Ricci tensor and scalar
curvature of a standard static space-time.

We recall the definition of a warped product of two
pseudo-Riemannian manifolds $(B,g_B)$ and $(F,g_F)$ with a
smooth function $b \colon B \to (0,\infty)$ (see \cite{BEE,ON}).
Suppose that $(B,g_B)$ and $(F,g_F)$ are pseudo-Riemannian
manifolds and also suppose that $b \colon B \to (0,\infty)$ is a
smooth function. Then the (singly) warped product, $B \times
{}_b F$ is the product manifold $B \times F$ equipped with the
metric tensor $g=g_B \oplus b^{2}g_F$ defined by $$
g=\pi^{\ast}(g_B) \oplus (b \circ \pi)^2 \sigma^{\ast}(g_F)$$
where $\pi \colon B \times F \to B$ and $\sigma \colon B \times
F \to F$ are the usual projection maps and ${}^\ast$~denotes the
pull-back operator on tensors. Here, $(B,g_F)$ is called as the
base manifold and $(F,g_F)$ is called as the fiber manifold and
also $b$ is called as the warping function. There are also
different generalizations of warped products such as warped
products with more than one fiber manifold, called multiply
warped products (see \cite{MWP}) or warped products with two
warping functions acting symmetrically on the fiber and base
manifolds, called doubly warped products (see \cite{DWP}).
Finally, a warped product is said to be a twisted product if the
warping function defined on the product of the base and fiber
manifolds (see \cite{TWP}).

Basically, a standard static space-time can be considered as a
Lorentzian warped product where the warping function is defined
on a Riemannian manifold and acting on the negative definite
metric on an open interval of real numbers. More explicitly, a
standard static space-time, ${}_f(a,b) \times F$ is a Lorentzian
warped product furnished with the metric $g=-f^2{\rm d}t^2
\oplus g_F,$ where $(F,g_F)$ is a Riemannian manifold, $f \colon
F \to (0,\infty)$ is smooth, and $-\infty \leq a < b \leq
\infty.$ This class of space-times has been previously
considered by many authors. Now, we give a summary of some major
work about standard static space-times. In \cite{ON}, it was
shown that any static space-time is locally isometric to a
standard static space-time. Kobayashi and Obata \cite{KO} stated
the geodesic equation for this class of space-times and the
causal structure and geodesic completeness was considered in
\cite{AD}, where sufficient conditions on the warping function
for nonspacelike geodesic completeness of the standard static
space-time was obtained (see also \cite{RASM} and \cite{DWP}).
In \cite{AD1}, conditions are found which guarantee that
standard static space-times either satisfy or else fail to
satisfy certain curvature conditions from general relativity.
The existence of geodesics in standard static space-times have
been studied by several authors. S\'{a}nchez \cite{mS2} gives a
good overview of geodesic connectedness in semi-Riemannian
manifolds, including a discussion for standard static
space-times. In \cite{GES}, geodesic structure of standard
static space-times is studied and conditions are found on the
warping function to imply non-returning and pseudo-convex
geodesic systems on a standard static space-time.

The Minkowski space-time and the Einstein static universe are
two most famous examples of standard static space-times (see
\cite{BEE,HE}) which is $\mathbb R \times \mathbb S^3$ equipped
with the metric
$$g=-{\rm d}t^2+({\rm d}r^2+\sin^2 r {\rm d}\theta^2+
\sin^2 r \sin^2 \theta {\rm d} \phi^2)$$ where $\mathbb S^3$ is
the usual 3-dimensional Euclidean sphere and the warping
function $f \equiv 1.$ Another well-known example is the
universal covering space of anti-de Sitter space-time, a
standard static space-time of the form ${}_f \mathbb R \times
\mathbb H^3$ where $\mathbb H^3$ is the 3-dimensional hyperbolic
space with constant negative sectional curvature and the warping
function $f \colon \mathbb H^3 \to (0,\infty)$ defined as
$f(r,\theta,\phi)=\cosh r$  (see \cite{BEE,HE}). As a final
example, we can also mention the Exterior Schwarzschild
space-time (see \cite{BEE,HE}), a standard static space-time of
the form ${}_f \mathbb R \times (2m,\infty) \times \mathbb S^2,$
where $\mathbb S^2$ is the 2-dimensional Euclidean sphere, the
warping function $f \colon (2m,\infty) \times \mathbb S^2 \to
(0,\infty)$ is given by $f(r,\theta,\phi)=\sqrt{1-2m/r},$ $r>2m$
and the line element on $(2m,\infty) \times \mathbb S^2$ is
$${\rm d}s^2=\bigl(1-\frac{2m}{r} \bigl)^{-1} {\rm
d}r^2+r^2({\rm d} \theta^2+\sin^2 \theta {\rm d} \phi^2).$$

In literature, it is not known in full generality that whether
there exists a nontrivial warping function $f$ for which the
warped product pseudo-Riemannian (or Riemannian) manifold $_fB
\times F$ is Einstein for given pseudo-Riemannian (or
Riemannian) manifolds $(F,g_F)$ and $(B,g_B)$.  In fact, the
answer of this question depends on whether there exists a
nontrivial common solution $f$ to some differential equations on
$(F,g_F)$. (See page 267 of \cite{B}). This problem was
considered especially for Einstein Riemannian warped products
with compact base and some partial answers were also provided
(see \cite{GA,KD,KK, KC}). In \cite{KK}, it is proven that an
Einstein Riemannian warped product with a non-positive scalar
curvature and compact base is just a trivial Riemannian product.
Constant scalar curvature of warped products was studied in
\cite{CDM,DD,EJK,EJKS} when the base is compact and of
generalized Robertson-Walker space-times in \cite{EJK}.
Furthermore, partial results for warped products with
non-compact base were obtained in \cite{BD} and \cite{CB}. The
physical motivation of existence of a positive scalar curvature
comes from the positive mass problem. More explicitly, in
general relativity the positive mass problem is closely related
to the existence of a positive scalar curvature (see \cite{ZX}).
As a more general related reference, one can consider \cite{KJ}
to see a survey on scalar curvature of Riemannian manifolds.

The problem of existence of a warping function which makes the
warped product Einstein was already studied for special cases
such as generalized Robertson-Walker space-times and a table
given summarizing different cases of Einstein Ricci tensor of a
generalized Robertson-Walker when the Ricci tensor of the fiber
is Einstein in \cite{ARS} (see also references therein). In this
paper, we consider this problem for standard static space-times.
In fact, we essentially investigate the conditions for $(F,g_F)$
so that there exists no nontrivial function $f$ on $F$
guaranteing that the standard static space-time $_f(a,b) \times
F$ is Einstein. Although the results in this paper remain valid
for the Riemannian setting, that is, when $(B,g_B)$ is the
Euclidean interval $((a,b),{\rm d}t^2)$, we prefer to state
these results in Lorentzian setting since there are certain
standard static space-times in Lorentzian geometry, such as,
Einstein static universe, Schwarzschild exterior space-time and
(universal) anti-de Sitter space-time, which are of interest in
Relativity Theory and appear as examples to the results of this
paper.  See {\it Remark \ref{rem3}} for the summary of the
results of this paper from the view of this Introduction.

Einstein Ricci tensor and constant scalar curvature of standard
static space-times with perfect fluid were already considered in
(see \cite{KO,MM}). Note that a matter is called a perfect fluid
if the energy-momentum tensor ${\rm T}$ has the form
$${\rm T}=(\mu+p)W \otimes W+pg,$$ where $W$ is a 1-form with
$g(W,W)=-1$ and $\mu$ and $p$ are called the energy density and
the pressure, respectively. In \cite{KO}, it is shown that a
standard static space-time ${}_f(a,b) \times F$ has a perfect
fluid if $${\rm Ric}_F-\frac{\tau_F}{s}g_F=\frac{1}{f}
\Bigl({\rm H}_F^f -\frac{\Delta_F(f)}{s}g_F \Bigl),$$ where
${\rm Ric}$ is the Ricci tensor, ${\rm H}$ is the Hessian form
and $\tau_F$ is the scalar curvature on $(F,g_F)$ and also ${\rm
dim}(F)=s$ (see also \cite{HE,ON,MM}). Moreover, in \cite{KO2},
the conformal tensor on standard static space-times with perfect
fluid is studied and it is shown that a standard static
space-time with perfect fluid is conformally flat if and only if
its fiber is Einstein and hence of constant curvature.

In this paper, we will consider arbitrary standard static
space-times, i.e., we will not assume the existence of a perfect
fluid, our results can be considered as extensions of the
results in \cite{KO,KO2,MM} where standard static space-times
with perfect fluid were considered. Duggal studied the scalar
curvature of 4-dimensional triple Lorentzian products of the
form $L \times B \times {}_fF$ and obtained explicit solutions
for the warping function $f$ to have a constant scalar curvature
for this class of products (see \cite{DK}). We also discuss
conditions on the warping function or on the fiber of a standard
static space-time to have a constant scalar curvature on the
space-time. Especially, we show that the Einstein static
universe cannot be generalized to a standard static space-time
$_f\mathbb{R}\times \mathbb{S}^s$ modelling a ``static universe"
with a nonconstant warping function $f$ on $(S^s,d\sigma^2)$,
that is, to a standard static space-time with constant scalar
curvature in {\it Theorem \ref{thm1}}.

\section{Preliminaries}

In this section, we give the formal definition of a standard
static space-time and state some curvature formulas (see
\cite{BEE,ON}).

\begin{dfn} \label{dfna} Let $(F,g_F)$ be an $s$-dimensional
Riemannian manifold and $f \colon F \to (0,\infty)$ be a smooth
function. Then $n(=s+1)$-dimensional product manifold
$(a,b)\times F$ furnished with the metric tensor $g=-f^2{\rm
d}t^2 \oplus g_F$ is called a standard static space-time and is
denoted by $_f (a,b) \times F$, where ${\rm d}t^2$ is the
Euclidean metric tensor on $(a,b)$ and $-\infty \leq a < b \leq
\infty.$
\end{dfn}

Throughout the paper the fiber $(F,g_F)$ of a standard
static space-time of the form $_f(a,b) \times F$ is always
assumed to be connected. Now we state some curvature
formulas for standard static space-times to be used later
in the proofs. Note that, since $_f (a,b) \times F$ and
$F\times_f(a,b)$ are isometric, the curvature formulas
below can easily be obtained from the well known curvature
formulas for warped product metric tensors by making
suitable substitutions (see for example,
\cite{BEE,B,ON,DWP}).

Here, we use sign convention for the Laplacian in \cite{ON},
i.e., defined by or $\Delta={\rm tr}({\rm H}),$ (see page 86 of
\cite{ON}) where ${\rm H}$ is the Hessian form (see page 86 of
\cite{ON}) and ${\rm tr}$ denotes for the trace, or
equivalently, $\Delta={\rm div}({\rm grad}),$ where ${\rm div}$
is the divergence and ${\rm grad}$ is the gradient (see page 85
of \cite{ON}).

\begin{prp} \label{prpa} Let $_f(a,b)\times F$ be a standard
static space-time. If $\tau$ and $\tau_F$ denote the scalar
curvatures of the space-time and fiber, respectively, then
$$\tau=\tau_F-2\frac{\Delta_F(f)}{f},$$
where $\Delta_F$ denotes the Laplace operator on $F.$
\end{prp}

\begin{rem} \label{rema} Note that, in the above Proposition,
since $\tau (t,p)$ is independent of $t\in (a,b)$, $\tau$ is a
lift of a unique function $\tilde\tau$ on $F$ to $(a,b)\times F$.
For brevity in expressions, sometimes we abuse the notation and
write $\tau$ instead of $\tilde\tau$ on $F$ to avoid taking the
lifts of functions on $F$ to $(a,b)\times F$.
\end{rem}

\begin{prp} \label{prpb} Let $_f(a,b)\times F$ be a standard
static space-time and also let $V$ and $W$ be vector fields on
$F.$ If ${\rm Ric}$ and ${\rm Ric}_F$ denote the Ricci tensors
of the space-time and the fiber, respectively, then
$${\rm \Ric}(\frac{\partial}{\partial t}+V, \frac{\partial}{\partial
t}+W)={\rm Ric}_F(V,W) + f \Delta_F(f)- \frac{1}{f}{\rm
H}^f_F(V,W),$$ where ${\rm H}^f_F$ is the Hessian form of $f$ on
$(F,g_F).$
\end{prp}

\section{Standard Static Space-times of Constant Scalar
Curvature}

In this Section, we essentially investigate geometric and
topological conditions on a Riemannian manifold $(F,g_F)$
which yield the constancy of the warping function $f$ on
$F$ of a constant scalar curvature standard static
space-time $_f(a,b)\times F$.

\begin{thm} \label{thm1} Let $_f(a,b)\times F$ be a standard
static space-time where $s\geq 2$. Assume that $(F,g_F)$ is
compact and the scalar curvature $\tau$ of the space-time is
constant. Then, $(F,g_F)$ is of constant scalar curvature
$\tau_F$ if and only if $f$ is constant on $F$. In either case,
$\tau=\tau_F.$
\end{thm}

\begin{proof} First note that, $\int_F(\tau_F-\tau)f=0$ by
{\it Proposition \ref{prpa}} and page 104 of \cite{NPR}. Thus,
since $f>0$ on $M$, $\tau_F(p_0)=\tau$ at some $p_0\in F$. Now,
if $\tau_F$ is constant on $F$ then $\tau_F(p)=\tau$ for all
$p\in F$, and it follows from {\it Proposition \ref{prpa}} that
$\Delta_F(f)=0$ on $(F,g_F)$. Thus, $f$ is constant on $F$.
Conversely, if $f$ is constant on $F$ then, by {\it Proposition
\ref{prpa}}, $\tau_F(p)=\tau$ for all $p\in F$, and hence
$\tau_F$ is constant on $F$.
\end{proof}

Note that, one of the ingredients of the physical concept
of a ``static universe" (see \cite{F}) is the constancy of
the scalar curvature of the space-time modelling a ``static
universe", which in fact, corresponds to the constancy of
the trace of the stress-energy tensor of the space-time via
Einstein equation.  Note that, the Einstein static universe
$\mathbb{R}\times\mathbb{S}^s$ is of constant scalar
curvature, where $\mathbb{S}^s=(S^s,d\sigma^2)$ is the unit
Euclidean sphere. (See page 189 of \cite{BEE}).  Thus, by
{\it Theorem \ref{thm1}}, the Einstein static universe
cannot be generalized to a standard static space-time
$_f\mathbb{R}\times \mathbb{S}^s$ modelling a ``static
universe" with a nonconstant warping function $f$ on
$(S^s,d\sigma^2)$, that is, to a standard static space-time
with constant scalar curvature.

It is noticed that one cannot obtain a non-trivial standard
static space-time of a constant scalar curvature when it has a
compact fiber of a constant scalar curvature. Thus we should
focus on standard static space-times with compact fibers of
nonconstant scalar curvatures. In \cite{EJKS}, a similar problem
was considered on a wider class of warped products (see also
\cite{EJK}). In order to make use of \cite{EJKS}, we introduce a
linear operator ${\rm L} \colon {\rm H}^{1,2}(F) \to {\rm
H}^{1,2}(F)$ on a compact Riemannian manifold $(F,g_F)$ defined
by $${\rm L}(v)=-\Delta_F(v)+\frac{\tau_F(q)}{2}v,$$ where $v$
in the Sobolev space ${\rm H}^{1,2}(F).$ Then we are ready to
state the following result.

\begin{thm} \label{thm1a} Let $(F,g_F)$ be a compact Riemannian
manifold with variable scalar curvature $\tau_F \colon F \to
\mathbb R$ where $s\geq 2.$ Then there exists a smooth function
$f \colon F \to (0,\infty)$ such that the corresponding standard
static space-time $_f(a,b) \times F$ is of constant scalar
curvature $\tau$.
\end{thm}

\begin{proof} From {\it Proposition \ref{prpa}}, like in \cite{DD}
and \cite{EJKS}, we look for $\tau \in \mathbb R$ and $f \in
\mathcal C^\infty(F)$ such that ${\rm L}f=\tau f.$ It is well
known that this type of eigenvalue problem has only one
eigenvalue, $\lambda_1(\tau_F),$ (which is simple) such that the
corresponding eigenfunction is strictly positive (see
\cite{NPR}). Note the centrality of the compactness of $F.$
\end{proof}

\begin{rem} The constant scalar curvature of the standard
static space-time $\tau=2 \lambda_1$ where $\lambda_1$ is
the first eigenvalue of the operator ${\rm L}$ on ${\rm
H}^{1,2}(F).$
\end{rem}

Thus it is possible to produce a non-trivial standard static
space-time of constant scalar curvature when the fiber is
compact and has a nonconstant scalar curvature. Now we turn our
attention to the complete case (not necessarily compact),
roughly speaking, under suitable hypothesis for the curvature of
the fiber, we will give a necessary condition for constant
scalar curvature in a standard static space-time.

\begin{thm} \label{nccase-s} Let $(F,g_F)$ be a complete
manifold without boundary where $s \geq 2$ . Suppose the Ricci
curvature of $F$ is non-negative, and suppose $\Delta_F(\tau_F)
\leq 0$ and also $\|\nabla_F(\tau_F) \|=o(r(x)),$ where $r(x)$
is the distance from $x$ to some fixed point $p \in F.$ If
${}_f(a,b) \times F$ is a standard static space-time with
constant scalar curvature $\tau,$ then $\tau \leq \inf_F
(\tau_F).$
\end{thm}

\begin{proof} By contradiction, suppose that $\tau>\inf_F(\tau_F).$
{\it Proposition \ref{prpa}} implies that $$\Delta_F(f)-
\bigl(\frac{\tau_F}{2}-\frac{\tau}{2} \bigl)f=0,$$ with $\tau$
constant and $f$ is positive. Let $q=(\tau_F-\tau)/2,$ since
$\tau>\inf_F(\tau_F)$ there results $\inf_F(q)<0.$ thus by
Corollary 1.1 of \cite{LY}, we obtain a contradiction.
\end{proof}

\begin{rem} In the previous theorem, we may also require
$(F,g_F)$ be compact and in this case, by the variational
structure of the principal eigenvalue of the operator
$-\Delta_F+\tau_F/2,$ namely $$\frac{\tau}{2}=\inf_{{\rm
H}^{1,2}(F)}\frac{\int_F \|\nabla_F(u)\|^2+
\frac{\tau_F}{2}u^2}{\int_F u^2},$$ there results
$$\frac{\tau}{2} \geq \inf_{{\rm H}^{1,2}(F)} \frac{\int_F
\|\nabla_F(u)\|^2+ \frac{\inf_F(\tau_F)}{2}u^2}{\int_F u^2}=
\frac{\inf_F(\tau_F)}{2}.$$ Here, notice that $\inf_{{\rm
H}^{1,2}(F)}\frac{\int_F \|\nabla_F(u)\|^2+
\frac{\inf_F(\tau_F)}{2}u^2}{\int_F u^2}$ is the principal
eigenvalue of $-\Delta_F+\inf_F(\tau_F)/2$ on the compact
manifold $F.$ So, $\tau \geq \inf_F(\tau_F).$ But by the above
theorem $\tau \leq \inf_F(\tau_F),$ thus $\tau=\inf_F (\tau_F)
\geq 0,$ since the latter inequality holds because of the
non-negative Ricci curvature of $F.$
\end{rem}

We now state a simple result for 2-dimensional standard
static space-times with constant scalar curvatures. Note
that if $M=_f(a,b) \times (c,d)$ is a 2-dimensional
standard static space-time with the metric tensor
$g=-f^2{\rm d}t^2+{\rm d}x^2,$ then its scalar curvature
$\tau$ is given by $\tau=-2f^{\prime \prime}/f,$ where
$-\infty \leq a < b \leq \infty$ and $-\infty \leq c < d
\leq \infty$ also $f \colon (c,d) \to (0,\infty)$ is
smooth.

\begin{prp} Let $_f(a,b) \times (c,d)$ be a 2-dimensional
standard static space-time. Then the scalar curvature
$\tau$ is constant if $f$ satisfies one of the followings
\begin{enumerate}
\item $\displaystyle{f(x)=c_1x+c_2,}$ for some
$c_1,c_2 \in \mathbb R$ when $\tau=0,$
\item $\displaystyle{f(x)=c_1\exp(\sqrt{-2\tau}x)+
c_2\exp(-\sqrt{-2\tau}x),}$ for some $c_1,c_2 \in \mathbb
R$ when $\tau<0,$
\item $\displaystyle{f(x)=c_1\cos(\sqrt{2\tau}x)+
c_2\sin(\sqrt{2\tau}x),}$ for some $c_1,c_2 \in \mathbb R$
when $\tau>0.$
\end{enumerate}
\end{prp}

One can compare the previous result with the characterization of
constant Gauss curvature revolution of surfaces embedded in
$\mathbb R^3$ (see the examples in page 66-67 of \cite{KW},
chapter 3 of \cite{PP} and page 169 of \cite{doC}).

\section{Einstein Standard Static Space-times}

We now concentrate on the Ricci tensor of standard static
space-times. More precisely, we will try to determine
conditions on the warping function of a standard static
space-time so that the space-time becomes Einstein or
Ricci-flat when the fiber is Einstein or Ricci-flat.

Recall that an arbitrary $n$-dimensional pseudo-Riemannian
manifold $(M,g)$ is said to be Einstein with $\lambda$ if there
exists a smooth map $\lambda \colon M \to \mathbb R$ such that
$\Ric=\lambda g.$ Furthermore, if $(M,g)$ is Einstein with
$\lambda$ and $\dim(M)=n \geq 3,$ then $\lambda$ is constant and
$\lambda=\tau/n,$ where $\tau$ is the (constant) scalar
curvature of $(M,g).$ Also note that for a 2-dimensional
Einstein manifold $(M,g)$ with $\lambda,$ one cannot necessarily
conclude the constancy of $\lambda$ (see \cite{ON}).

\begin{prp} \label{prp1}  Let $(F,g_F)$ be an $s$-dimensional
Riemannian manifold with scalar curvature $\tau_F$ where $s\geq
2$ and let the standard static space-time $_f(a,b)\times F$ be
Einstein with constant scalar curvature $\tau$. Then,
\begin{enumerate}
\item[a)] $\Delta_F(f)=-\frac{\tau}{n}f$ and
${\rm Ric}_F=\frac{1}{f}{\rm H}_F^f+\frac{\tau}{n}g_F$ on
$(F,g_F)$, where $n=s+1$.
\item[b)] $(s+1)\tau_F=(s-1)\tau$.
\end{enumerate}
\end{prp}

\begin{proof} By assumption, since
${\rm Ric}=\frac{\tau}{n}(-f^2dt^2\oplus g_F)$, we obtain from
{\it Proposition \ref{prpb}} that,
$${\rm Ric}_F(V,W)+f\Delta_F(f)-\frac{1}{f}{\rm H}_F^f(V,W)=-
\frac{\tau}{n}f^2+\frac{\tau}{n}g_F(V,W)$$ for all $V,W\in \Gamma
TF$.
\begin{enumerate}
\item[{\it a)}] If we set $V=0=W$ in the above expression, we obtain
$\Delta_F(f)=-\frac{\tau}{n}f$ and hence, it also follows
that ${\rm Ric}_F(V,W)=\frac{1}{f}{\rm
H}_F^f(V,W)+\frac{\tau}{n}g_F(V,W)$ for all $V,W\in\Gamma
TF$.
\item[{\it b)}] Note that, if we take the trace of the equation
${\rm Ric}_F=\frac{1}{f}{\rm H}_F^f+\frac{\tau}{n}g_F$ with
respect to $g_F$ on $F$, we obtain
$\tau_F=\frac{1}{f}\Delta_F(f)+\frac{s}{n}\tau$, and hence,
it follows from {\it (a)} that $(s+1)\tau_F=(s-1)\tau$.
\end{enumerate}
\end{proof}

\begin{rem} Let $(F,g_F)$ be a Riemannian manifold of
scalar curvature $\tau_F$ where $s\geq 2$ and let the
standard static space-time $_f(a,b)\times F$ be Einstein
with constant scalar curvature $\tau$.  Then, by {\it
Theorem \ref{thm1}} and {\it Proposition \ref{prp1}}, we
conclude the following:
\begin{enumerate}
\item[a)] $(F,g_F)$ is of constant scalar curvature $\tau_F$.
\item[b)] If $(F,g_F)$ is compact then $\tau=0,$ \,
$\tau_F=0$ and $f$ is constant on $F$.
\end{enumerate}
\end{rem}

\begin{rem} \label{rem1} Note that, in {\it Proposition
\ref{prp1}}, if we further assume that $(F,g_F)$ is
Einstein with scalar curvature $\tau_F$ then, since ${\rm
Ric}_F=\frac{\tau_F}{s}g_F$ on $(F,g_F)$, we obtain by
using {\it Proposition \ref{prp1}} that ${\rm
H}_F^f=(\frac{\tau_F}{s}-\frac{\tau}{n})fg_F=-\frac{\tau_F}{s(s-1)}g_F
=-\frac{\tau}{s(s+1)}g_F$ on $(F,g_F)$.  (Note that, by
{\it Proposition \ref{prp1}-b}, $\tau_F$ is also constant
when $dim\,F=s=2$).
\end{rem}

\begin{thm} \label{thm2} Let $(F,g_F)$ be a complete
Riemannian manifold with nonnegative Ricci curvature where
$s\geq 2$. If the standard static space-time $_f(a,b)\times
F$ is Ricci flat then $f$ is constant on $F$.
\end{thm}

\begin{proof} Since $_f(a,b)\times F$ is Einstein with scalar
curvature $\tau = 0$, it follows from {\it Proposition
\ref{prp1}} that $\Delta_F(f)=0$ on $(F,g_F)$.  Since $f$
is positive on $F$, $f$ is constant by {\it Corollary 1} of
\cite{Y}.
\end{proof}

Note that, by {\it Theorem \ref{thm2}}, if a standard
static space-time $_f(a,b)\times F$ is Ricci flat with a
nonconstant warping function $f$ on $F$ then $(F,g_F)$ is
either incomplete or not of nonnegative Ricci curvature (or
both).  Hence, for the Schwarzschild exterior space-time
$_f\mathbb{R}\times F$ (see page 367 of \cite{ON}), we
conclude that $(F,g_F)$ is either incomplete or not of
nonnegative Ricci curvature, where $N=(2m,\infty)\times
S^2$ and $g_F=\bigl(1-\frac{2m}{r} \bigl)^{-1}
dr^2+d\sigma^2$. Indeed, for the Schwarzschild exterior
space-time, $(F,g_F)$ is both incomplete and not of
nonnegative Ricci curvature.

\begin{rem} Note that the converse of {\it Theorem
\ref{thm2}} is not true in general.  For example, the
$n(\geq 3)$-dimensional Einstein static universe is an
counterexample to this case.
\end{rem}

\begin{thm} \label{thm3} Let $(F,g_F)$ be a complete
Einstein $s(\geq 2)$-dimensional Riemannian
manifold with scalar curvature $\tau_F$. If the standard
static space-time $_f(a,b)\times F$ is Einstein with scalar
curvature $\tau$ then $\tau\leq 0$. Furthermore,
\begin{enumerate}
\item[a)] if $\tau =0$ then $f$ is constant on $F$,
\item[b)] if $\tau <0$ then $f$ is nonconstant on $F$ and
$(F,g_F)$ is a warped product of the Euclidean line and a complete
Riemannian manifold with warping function $\psi$ on the real line
satisfying the equation
$\frac{d^2\psi}{dt^2}+\frac{\tau}{s(s+1)}\psi=0$, $\psi >0$.
\end{enumerate}
\end{thm}

\begin{proof} Suppose that $\tau >0$.  Then by {\it Proposition
\ref{prp1}}, $\tau_F=\frac{s-1}{s+1}\tau >0$ and it follows
from Myers theorem that $(F,g_F)$ is compact.  Hence, by
{\it Theorem \ref{thm1}}, $f$ is constant on $F$.  But this
conflicts with $\Delta_F(f)=-\frac{\tau}{n}f$ since $f>0$
on $F$ (see {\it Proposition \ref{prp1}}).  Thus $\tau \leq
0$.  Furthermore,
\begin{enumerate}
\item[{\it a)}] if $\tau =0$ then, by {\it Proposition \ref{prp1}},
$\Delta_F(f)=0$ on $(F,g_F)$ and ${\rm
Ric}_F=\frac{\tau_F}{s}g_F=0$ since $\tau_F=0$.  Thus,
since $f$ is positive on $F$, $f$ is constant by {\it
Corollary 1} of \cite{Y}.
\item[{\it b)}] if $\tau <0$ then, by {\it Proposition
\ref{prp1}} and {\it Remark \ref{rem1}}, $f$ is nonconstant
and ${\rm H}_F^f=-\frac{\tau_F}{s(s-1)}g_F$ on $(F,g_F)$,
where $\tau_F=\frac{s-1}{s+1}\tau <0$.  Hence, it follows
from {\it Corollary E} of \cite{KA} that $(F,g_F)$ is a
warped product of the Euclidean line and a complete
Riemannian manifold with warping function $\psi$ on the
real line satisfying the equation
$\frac{d^2\psi}{dt^2}+\frac{\tau}{s(s+1)}\psi=0$, $\psi
>0$.
\end{enumerate}
\end{proof}

Note that, the (universal) anti-de Sitter space-time
$_f\mathbb{R}\times\mathbb{H}^s$ of constant sectional
curvature $-1$ is an example to {\it Theorem \ref{thm3}-b}
for $s\geq 2$. (See page 183 of \cite{BEE}).  Indeed, note
that the Riemannian hyperbolic space $\mathbb{H}^s$ can be
written as a warped product of the Euclidean line and the
Euclidean space with warping function $\psi =e^{\pm t}$ on
the real line.

As we did in {\it Section 3} we will consider the same type of
problems on a standard static space-time by using \cite{LY} when
the fiber is complete without boundary.

\begin{thm} \label{nccase-e1} Let $(F,g_F)$ be a complete
manifold without boundary where $s \geq 2.$ Suppose the Ricci
curvature of $F$ is non-negative, and suppose $\Delta_F(\tau_F)
\leq 0$ and also $\|\nabla_F(\tau_F) \|=o(r(x)),$ where $r(x)$
is the distance from $x$ to some fixed point $p \in F.$ If
${}_f(a,b) \times F$ is an Einstein standard static space-time
then $\tau \leq \frac{s+1}{s}\inf_F(\tau_F).$
\end{thm}

\begin{proof} First note that the space-time is Einstein this
means that ${\rm \Ric}=\lambda g,$ where $(s+1)\lambda=\tau$ and
$\tau$ is constant. Considering the trace in $(F,g_F)$ in {\it
Proposition \ref{prp1}(a)}, there results a positive solution
$f$ for the Schr\"odinger equation $-\Delta_F(f)+q_E(f),$ where
$q_E=\tau_F-\frac{s}{s+1} \tau.$ On the other hand, $q_E$
verifies $\Delta_F(q_E) \leq 0$ and $\|\nabla_F(q_E)\|=o(r(x))$
with $r(x)$ as in the hypothesis. So by Corollary 1.1 of
\cite{LY}, $0 \geq \inf_F(q_E)=\inf_F-\frac{s}{s+1}\tau,$ or
equivalently $\tau \leq \frac{s+1}{s} \inf_F(\tau_F).$
\end{proof}

\begin{thm} \label{nccase-e2} Let $(F,g_F)$ be a complete
manifold without boundary where $s \geq 2.$ Suppose the Ricci
curvature of $F$ is non-negative. If ${}_f(a,b) \times F$ is an
Einstein standard static space-time then $\tau=\tau_F=0$ and
hence the space-time is Ricci-flat.
\end{thm}

\begin{proof} As again in the proof of the previous theorem,
note that the space-time is Einstein this means that ${\rm
\Ric}=\lambda g,$ where $(s+1)\lambda=\tau$ and $\tau$ is
constant. By  {\it Proposition \ref{prp1}(a)},
$\tau_F=\frac{s-1}{s+1} \tau$ is constant. Thus the hypothesis
of {\it Theorem \ref{nccase-e1}} are satisfied and this leads
$\tau \leq \frac{s+1}{s} \inf_F(\tau_F)=\frac{s-1}{s} \tau.$ On
the other hand, as Ricci curvature of $F$ is non-negative,
$\tau=\frac{s+1}{s-1}\tau_F \geq 0.$ So, $\tau=0.$
\end{proof}

\begin{cor} \label{nccase-e3} Let $(F,g_F)$ be a complete
manifold without boundary where $s \geq 2.$ Suppose the Ricci
curvature of $F$ is non-negative. If ${}_f(a,b) \times F$ is an
Einstein standard static space-time then $f$ is constant and
$(F,g_F)$ is Ricci-flat.
\end{cor}

\begin{proof} It just follows {\it Theorem \ref{thm2}}
and {\it Theorem \ref{nccase-e2}} and also {\it Proposition
\ref{prp1}(a)}.
\end{proof}

As a consequence, we obtain the result that follows:

\begin{thm} \label{nccase-e4} Let $(F,g_F)$ be a complete
manifold without boundary where $s \geq 2.$ If ${}_f(a,b) \times
F$ is an Einstein standard static space-time then either
$(F,g_F)$ is Ricci-flat or the Ricci curvature of $F$ cannot be
non-negative.
\end{thm}

In the remark below, we collect the results in this section from
viewpoint of the nonexistence of a nonconstant warping function
$f$ on a connected, complete Riemannian manifold $(F,g_F)$ for
which the standard static space-time $_f(a,b)\times F$ is
Einstein.

\begin{rem} \label{rem3} Let $(F,g_F)$ be a Riemannian
manifold and $_f(a,b)\times F$ be a standard static space-time
where $s\geq 2.$
\begin{enumerate}
\item[a)] If $(F,g_F)$ is compact and of constant scalar curvature
(or Einstein) then there exists no nonconstant function $f$ on
$F$ for which $_f(a,b)\times F$ is of constant scalar curvature
(or Einstein). (See {\it Theorem \ref{thm1}}).
\item[b)] If $(F,g_F)$ is complete and of nonnegative Ricci
curvature then there exists no nonconstant function $f$ on $F$
for which $_f(a,b)\times F$ is Ricci flat. (See {\it Theorem
\ref{thm2}}).
\item[c)] If $(F,g_F)$ is complete and Einstein then there exists
no nonconstant function $f$ on $F$ for which $_f(a,b)\times F$
is Ricci flat. (See {\it Theorem \ref{thm3}}).
\item[d)] If $(F,g_F)$ is complete and of non-negative Ricci
curvature, then there exists no nonconstant function $f$ on $F$
for which $_f(a,b)\times F$ is Einstein. (See {\it Corollary
\ref{nccase-e3}}).
\end{enumerate}
\end{rem}

Here note that, the above Remark remains valid (as well as other
results in this section) in the Riemannian setting, that is,
when we take $(a,b)$ as an (Euclidean) interval with metric
tensor ${\rm d}t^2$.

We now consider 2-dimensional Einstein standard static
space-times. Let $M=_f(a,b) \times (c,d)$ be a 2-dimensional
standard static space-time with the metric tensor $g=-f^2{\rm
d}t^2+{\rm d}x^2.$ If $\frac{\partial}{\partial t} \in \mathfrak
X(a,b)$ and $\frac{\partial}{\partial x} \in \mathfrak X(c,d),$
then
$${\rm Ric}(\frac{\partial}{\partial t}+
\frac{\partial}{\partial x})=ff^{\prime \prime}-\frac{f^{
\prime \prime}}{f},$$ where $-\infty \leq a < b \leq
\infty$ and $-\infty \leq c < d \leq \infty$ also $f \colon
(c,d) \to (0,\infty)$ is smooth. Thus we can easily state
the following result.

\begin{prp} Let $_f(a,b) \times (c,d)$ be a 2-dimensional
standard static space-time. Then
\begin{enumerate}
\item the space-time is Einstein with $\lambda$ if and only
if $\displaystyle{f^{\prime \prime}=-\lambda f},$
\item the space-time is Ricci-flat if and only
if $\displaystyle{f(x)=c_1x+c_2}$ on $(c,d),$ for some
$c_1,c_2 \in \mathbb R.$
\end{enumerate}
\end{prp}

In the previous result, since we cannot conclude the constancy
of $\lambda,$ it is impossible for us to obtain explicit
solutions. However, we can only say that $\lambda \colon (a,b)
\times (c,d) \to \mathbb R$ depends only on the second variable,
i.e., for each $x \in (c,d),$ we have $\lambda(t_1,x)=
\lambda(t_2,x)$ for any $t_1, t_2 \in (a,b).$

\vspace{0.5 cm}

{\sc Authors'} {\sc Addresses:}

\vspace{0.5 cm}

\begin{minipage}[b]{14 cm}
\parbox[t]{8.0 cm}{{\bf Fernando Dobarro}\\
Dipartimento di Scienze Matematiche \\
Universit\`{a} degli Studi di Trieste \\
Via Valerio 12 I-34127, Trieste \\
Italy \\ {\it e-mail:} dobarro@mathsun1.univ.trieste.it}
\hfill \hfill \hfill \hfill
\parbox[t]{6.0 cm}{{\bf B\"ulent \"Unal}\\
Department of Mathematics \\
Atilim University \\
Incek 06836, Ankara \\
Turkey \\
{\it e-mail:} bulentunal@mail.com}
\end{minipage}

\end{document}